# Restarting the Numerical Flow Iteration through low rank tensor approximations

Rostislav-Paul Wilhelm[*,†]
Katharina Kormann[‡]

October 10, 2025

**Abstract**

The numerical flow iteration method has recently been proposed as a memory-slim solution method for the Vlasov–Poisson system. It stores the temporal evolution of the electric field and reconstructs the solution in each time step by following the characteristics backwards in time and reconstructing the solution from the initial distribution. If the number of time steps gets large, the computational cost of this reconstruction may get prohibitive. Given a representation of the intermediate solution, the time intervals over which the characteristic curves need to be solved backwards in time can be reduced by restarting the Numerical Flow Iteration (NuFI) after certain time intervals. In this paper, we propose an algorithm that reconstructs a low-rank representation of the solution at the restart times using a randomized singular value decomposition (RSVD) algorithm with lazy evaluation. The proposed algorithm reduces the computational complexity compared to the pure Numerical Flow Iteration from quadratic to linear in the number of times step while still keeping its memory complexity. On the other hand, our numerical results demonstrate that the methods preserves the property of the Numerical Flow Iteration of showing much less dissipation of filaments compared to the semi-Lagrangian method.

Keywords: Plasma Physics, Kinetic Theory, Vlasov equation, Numerical Methods, Low Rank Tensor approximation

## 1 Introduction

In many applications such as nuclear fusion devices or solar wind it is relevant to fundamentally understand how turbulence and instabilities in a plasma are triggered and how they evolve over time. Some phenomena can already be observed through a fluid-description but when plasma becomes rarefied or hot the velocity distribution of the charged particles can substantially deviate from equilibrium, which breaks the fundamental assumption of a fluid-descriptions. In this case one has to model the plasma evolution through kinetic theory, i.e., including the velocity distributions in the model, which results in the Vlasov equation[1]

$$\partial_t f^\alpha + v \cdot \nabla_x f^\alpha + F^\alpha \cdot \nabla_v f^\alpha = 0. \tag{1}$$

where $f^\alpha$ is the probability distribution of the species $\alpha$ in the up to six-dimensional phase-space and $F^\alpha$ are the self-induced and external electro-magnetic forces acting upon the species $\alpha$. In this work we will only consider the single-species, electro-static case as the complications occurring in data-compression, which we want to discuss here, are already sufficiently pronounced in this simplified setup.

[*]Centre for mathematical Plasma Astrophysics, KU Leuven, Belgium
[†]Institute for Applied and Computational Mathematics, RWTH Aachen University
[‡]Department of Mathematics, Ruhr University Bochum, Bochum

Hence, in the following we will drop the subscript $\alpha$ and restrict ourselves to the self-induced electric field force $F = \frac{q}{m}E$. In this case we can additionally assume $\frac{q}{m} = 1$ without loss of generality. This coupling yields the non-linear Vlasov–Poisson system. Finally, if not explicitly stated otherwise, we assume an uniform ion background and only consider electron dynamics, thus we only have to solve a single Vlasov equation coupled to the Poisson equation

$$-\Delta_x \varphi = \rho = 1 - \int_{\mathbb{R}^3} f \mathrm{d}v, \tag{2}$$

$$E = -\nabla_x \varphi. \tag{3}$$

As any hyperbolic conservation law, the solution to (1) remains constant along the characteristic curves given as

$$\begin{aligned} \tfrac{\mathrm{d}}{\mathrm{d}s}\hat{x}(s) &= -\hat{v}(s), & \hat{x}(t) &= x, \\ \tfrac{\mathrm{d}}{\mathrm{d}s}\hat{v}(s) &= -E\left(s, \hat{x}\left(s\right)\right), & \hat{v}(t) &= v. \end{aligned} \tag{4}$$

Therefore, the solution to (1) can be written as

$$f(t, x, v) = f_0(\Phi_t^0(x, v)), \tag{5}$$

where the *backward flow* is given by $s \mapsto \Phi_t^s(x, v) = (\hat{x}(s), \hat{v}(s))$.

Solving the Vlasov system poses a number of significant challenges: The solution lives in the six-dimensional phase-space inducing the curse of dimensionality for grid-based methods. Therefore simulations of the system are only possible through either low resolution or with a number of additional assumptions on the physics to reduce the dimension of phase space. Additionally, the non-linearity and lack of collisions in the Vlasov system pose further complications. The Vlasov system is notoriously known for the development of fine structures, called filaments, in the distribution function, which are hard to capture through conventional approaches but are relevant to capture the onset of kinetic instabilities and there is some evidence that they play a key role in the heat dissipation in a plasma.[2],[3],[4]

Solvers for the Vlasov equation can be divided into 3 general categories: Eulerian (or Discrete Velocity), Semi-Lagrangian and Lagrangian. The first two are grid-based, i. e., they store the solution $f$ to the Vlasov equation on a phase-space grid. While Eulerian approaches update the values on the grid through looking at the local rate of change, transforming the Vlasov equation into a set of ordinary differential equations (ODE) with as many ODEs as there are degrees of freedom in the discretization.[5],[6],[7],[8],[9],[10],[11] While Eulerian methods are common in the context of fluid dynamics, they are less commonly used to solve the Vlasov or related kinetic equations: the Eulerian viewpoint provides a “general purpose recipe” to derive a numerical solver for a partial differential equation (PDE), however, it usually does not take additional structural information into account. Because the Vlasov–Poisson – and Vlasov–Maxwell – equations are divergence-free, characteristics do not cross. The semi-Lagrangian method uses a phase-space grid but evolves the grid points backwards in time on small time intervals by the characteristic equations and then reconstructs the solution at the new time step by interpolation from the solution at the previous time step evaluated at the foot of the characteristic.[1] This improves the conservation properties of the scheme compared to a pure Eulerian approach, especially if a symplectic time-integration scheme is used in addition.[12] A common choice for interpolation are (cubic) B-Splines as they offer a good compromise between accuracy and performance[13],[14] but there are also a number of codes that use nodal interpolation method for the approximation.[15],[16],[17],[18],[19],[20],[21],[22],[23] However, grid-based approaches for the Vlasov equation face a number of issues: The curse of dimensionality makes maintaining a 6-dimensional grid prohibitively expensive, especially if adaptivity or non-rectangular domains are required. Communication overhead introduced through the high memory footprint limits parallel efficiency and thus also scalability, especially on modern accelerator hardware.[14],[24],[25],[26] Introducing sparse grids or adaptivity into the simulations are very challenging due to the high dimensions.[27],[28] Grids always introduce numerical diffusion, which in the case of the collisionless Vlasov equation reintroduces unphysical dissipation into the simulation and thus violates the inherent conservation laws of the system. While there are some works towards limiting this effect, it is inherent to any grid-based approach and therefore cannot be avoided altogether.[11],[16] Finally, handling of more complicated domain shapes, boundary conditions or the expanding velocity support is another major challenge for grid based approaches.[29],[30],[31]

Alternatively to a grid-based solver it is also possible to use particle-based – so-called Lagrangian – solvers. The most common approach is the *Particle-In-Cell* (PIC) method, though there are also some particle-based methods such as *smoothed particle hydrodynamics* (SPH) or the *Reproducing Kernel Hilbert Space Particle method* (RKHS-PM).[32],[33],[34],[35] Lagrangian methods represent the distribution function through a set of marker particles which carry function values or weights, which are moved along the phase flow associated to the velocity field of the Vlasov equation. As these trajectories of the particles are in fact the characteristics of the Vlasov equation, the respective function values (or weights) do not change over time (this is the same idea used by Semi-Lagrangian approaches). In PIC one additionally prescribes a grid in the physical domain to solve for and evaluate the electro-magnetic fields. To compute the charge density – the right-hand side of (2) – one can count the particles per cell or use higher-order shape functions to deposit the charge.[36],[37],[38],[39],[40] While PIC also suffers from the curse of dimensionality – even though to a smaller extend, the major advantage over grid-based approaches is that particles automatically adapt to dynamics. Furthermore, as data is unstructured in PIC, it is also easier to handle more complicated domains as well as boundary conditions and parallelization for distributed computing. However, particle-based representation suffer from numerical noise and thus one has to use higher resolutions – compared to grid-based solvers – to achieve the same target accuracy. In simulations with low density regions or where high accuracy is required the noise becomes an issue. Most common are explicit PIC codes,[41],[42],[43],[44],[45],[46] i. e., where the time integration is explicit. Due to the large difference in mass between electrons and ions in a plasma, plasma dynamics usually are a very multi-scale phenomena. Explicit codes often require prohibitively high resolutions to capture all relevant scales, hence there has also been some development of fully and semi-implicit PIC codes.[47],[48],[49],[50],[51],[52],[53]

A more detailed overview of numerical solvers for the Vlasov equation can be found in e.g. the reviews by Filbet and Sonnendrücker[54] or, more recently, Palmroth et al.[31]

[1]Technically the described scheme is a *backwards Semi-Lagrangian* (BSL) approach. If tracking grid points forwards in time, one obtains a *forward Semi-Lagrangian* (FSL) approach. For simplicity of notation in future we always refer to BSL when saying SL.

The major challenge in solving kinetic equations is the high dimensionality. In recent years there has been a push to address this through using *low rank compression algorithms*.[55] In particular, *dynamical low rank* (DLR), where the initial data is already prescribed in low rank format and then evolved by restricting the evolution to the manifold of low rank solutions, has gained some interest.[56],[57],[58],[59],[60],[61],[62],[63] For dynamics which are dominated by collisional processes, such as e.g. radiative transfer, this approach works very well as the collisions effectively limit the maximum rank.[62],[64],[65] However, as the Vlasov equation is collisionless fine structures develop over time, which have to be resolved to accurately capture the physics but naturally increase the rank of the solution. This leads to a build-up of significant errors when applying DLR to Vlasov and to an increased loss of structure, especially conservation properties, for long simulation times. While there are some attempts to limit the loss of conservation properties,[57],[63] it is inherent problem of DLR. Missing development of these fine scale structures can lead to missing relevant kinetic effects.

The difficulties with solving the Vlasov equation numerically can be explained by 3 interconnected key challenges: the high-dimensionality, (physically relevant) filamentation as well as turbulence and the multi-scale nature of the dynamics. So the challenge is to resolve high dimensional solutions while maintaining the structure and filamentation across multiple scales. Trying to directly discretize through resolving and storing the distribution functions is bound to require extremely high resolution, which in turn induce prohibitively high memory and computational cost, even with advanced sparse or adaptive meshing techniques. Hence with state-of-the-art algorithms the rate of physical discovery will always be directly bound to rate of hardware improvements.

We propose a change of paradigm: Instead of a direct discretization of $f$ we approximate the characteristics. This is possible in an efficient way through the novel *Numerical Flow Iteration* (NuFI) method as they can be reconstructed on-the-fly using operator-splitting. This only requires knowledge of the past electric fields, while maintaining the Hamiltonian structure of the Vlasov equation, so it not only reduces the memory cost by several orders of magnitude but also preserves the invariant of the analytic solution, i. e., conserves $L^p$-norms, kinetic entropy and total energy.[66] Previous publications discussed the efficient implementation of NuFI for single- and multi-species simulations in the electro-static limit with periodic boundary conditions[66],[67],[68] and non-periodic boundary conditions.[69] Due to the high degree of structure preservation by NuFI, the scheme is more accurate than other state-of-the-art methods, however, the paradigm change to evaluation on-the-fly comes at the cost of increasing the computational complexity from linear to quadratic in the number of time-steps for each grid point. This again makes long simulations with NuFI too computationally expensive.

To reduce the computational complexity back to linear the simplest approach would be to just restart NuFI periodically in time, i. e., every $n > 1$ time-steps store a snapshot of the distribution function, which then would act as the “new initial data”,[69] but this would reintroduce the previously discussed challenges, in particular, the high memory requirement. To avoid this we suggest to use low rank compression, however, only for compression of the solution data, not the evolution. The proposed combination of NuFI and low rank compression (NuFI-LR) would be comparable to a Semi- Lagrangian scheme with high fidelity subcycling via NuFI.

Note that recently another approach, closely related to NuFI, was proposed: the *Characteristic Mapping Method* (CMM). The idea is to directly discretize the characteristic flow map using a composition strategy preserving the semi-group structure. It can be shown that through composing maps which cover only short time spans, low resolutions are sufficient to resolve them accurately.[70],[71],[72] However, the characteristic flow map maps the 6-dimensional phase-space into itself, i. e., CMM has an even higher memory footprint than a classical Semi-Lagrangian approach, albeit more accurate and structure-preserving. Also as CMM does not use symplectic integrators, the conservation properties of the analytic solution are not preserved analytically as for NuFI. We are currently collaborating on merging the two methods, which will be reported in a forthcoming work.

In the following we will first briefly introduce the Numerical Flow Iteration in section 2. Then in section 3 we introduce the low rank tensor approximation of NuFI. Finally in section 4 we study the properties of the NuFI-LR method for three benchmark test cases and compare the results both the the NuFI method without restart and a standard Semi-Lagrangian solution.

## 2 The Numerical Flow Iteration for the Vlasov–Poisson equation

A detailed derivation and discussion of the algorithm can be found in previous work.[66],[67],[68] In this section we want to recap the important ideas and highlight the challenges, which we try to tackle in the following chapters.

The core idea of the Numerical Flow Iteration (NuFI) is to approximate the flow-map of the Vlasov–Poisson equation instead of directly solving for and storing the distribution functions $f$. This is possible as the evolution of the Vlasov–Poisson system, independent of the number of species, is completely described through the evolution of the electric field in time and the initial value through the representation (5) and the characteristic equations (4). To solve (4) numerically we use the Störmer–Verlet time-integration method: Starting from $\hat{x}_n^h = x$ and $\hat{v}_n^h = v$ at the time-step $t = t_n$ compute

$$\hat{v}_{i-1/2}^h = \hat{v}_i^h - \tfrac{\Delta t}{2} E(t_i, \hat{x}_i^{h,\alpha}), \tag{6}$$

$$\hat{x}_{i-1}^h = \hat{x}_i^h - \hat{v}_{i-1/2}^h, \tag{7}$$

$$\hat{v}_{i-1}^h = \hat{v}_{i-1/2}^h - \tfrac{\Delta t}{2} E(t_{i-1}, \hat{x}_{i-1}^h) \tag{8}$$

for $i = 0, ..., n$. Now

$$f(t, x, v) = f_0(\hat{x}_0^h, \hat{v}_0^h) + \mathcal{O}\left(\Delta t^2\right). \tag{9}$$

Note that to compute $\rho$ via (2) we can skip the first half-step (6) as it can be recast as a transformation of the respective integral with functional determinant 1. With this we can state the full NuFI algorithm (Algorithm 1). Note that we only write it out for $d = 1$ but extending to $d = 2, 3$ is trivial.[66]

**Algorithm 1** Solving the Vlasov–Poisson system in $d = 1$.

```
function NuFI(f_0^e, f_0^i, N_t, Δt, N_x, N_v)
    Allocate a array C for the coefficients of φ (N_t N_x floats).
    for n = 0, ..., N_t do
        for k = 0, ..., N_x do
            ρ_k^e, ρ_k^i = 0.
            for l = 0, ..., N_v do
                Evaluate f̃ = f(t_n, x_k, v_l) through (9) using (6) - (8) to solve (4).
                ρ_k += h_v f̃.
            end for
            ρ_k = 1 − ρ_k.
        end for
        Solve the Poisson's equation via FFT to obtain φ_0, ..., φ_{N_x} from ρ_0, ..., ρ_{N_x}.
        Interpolate φ_0, ..., φ_{N_x} and store the coefficients of φ_{Δt,h} in the B-Spline basis.
    end for
end function
```

It is no longer required to store the distribution functions explicitly and instead the approach only stores the coefficients of the lower dimensional electric potential. Therefore it has substantially lower memory requirement than other approaches which need to store the high dimensional distribution functions: The memory requirement of a “classical approach” needs to store

$$\mathcal{O}\left(N_x \cdot N_y \cdot N_z \cdot N_u \cdot N_v \cdot N_w\right) \tag{10}$$

degrees of freedom, while NuFI stores

$$\mathcal{O}\left(N_x \cdot N_y \cdot N_z \cdot N_t\right) \tag{11}$$

degrees of freedom.

Commonly $N_t \ll N_u \cdot N_v \cdot N_w$ holds. Thus we see NuFI's memory-advantage becomes more apparent the higher dimension we need to consider in our simulation and the more degrees of freedom in velocity are needed. In addition, as a symplectic time-integrator is employed it can be shown that NuFI analytically preserves the conservative structure of the Vlasov–Poisson system, which also leads to higher accuracy in simulations compared to other approaches using the same resolution in phase-space.[66]

However, these advantages are paid through the increased computational complexity of NuFI. While a classical kinetic solver has linear complexity in the number of time-steps, this is not the case for NuFI: While a classical approach needs to perform

$$\mathcal{O}\left(N_t \cdot N_x \cdot N_y \cdot N_z \cdot N_u \cdot N_v \cdot N_w\right) \tag{12}$$

operations for a simulation of $N_t$ time steps in total, NuFI has to perform

$$\mathcal{O}\left(N_t^2 \cdot N_x \cdot N_y \cdot N_z \cdot N_u \cdot N_v \cdot N_w\right) \tag{13}$$

operations. That is due to the method having to iterate all the way back to the initial data in each time-step for each evaluation of $f$, which is linear in the number of backwards iterations. As this is done for every time-step the overall complexity becomes quadratic in the number of total time steps.

Hence while simulations of higher resolution can be run due to the substantially lower memory requirement and they are both more accurate and can be executed more efficiently, these advantages are paid for with a higher computational complexity. This makes long simulation periods prohibitively expensive with pure NuFI.

# 3 Low rank tensor approximation of NuFI

As discussed in section 2, NuFI has improved accuracy and low memory consumption but this comes at cost of having to evaluate the entire characteristic map backwards until the initial data is reached, thus the computational complexity of NuFI is quadratic, not linear, in the total number of time-steps.

Considering that we are interested in kinetic instabilities and turbulence which can take substantial time periods to develop, it is crucial to reduce the computational complexity of NuFI back to linear. To this end we propose to restart NuFI periodically in time through e.g. storing an approximation to the distribution function $f_n$. Instead of always going back to $f_0$ one then goes back to $f_n$ for time-steps $m > n$. In a sense, this can be seen as a Semi-Lagrangian approach with sub-cycling done by NuFI without having to store intermediate results.[69] But to not reintroduce the memory bottleneck of having to store the full, high-dimensional distribution function again, it is crucial to use compression on the stored data as we will discuss in the second part of this section.

## 3.1 Compression by truncated singular value decomposition

Two-dimensional data that is represented on a tensor product of two one-dimensional grids of $m$ and $n$ points, respectively, can be represented as a matrix. The best rank-$k$ approximation of a matrix, is given by truncating its *singular value decomposition* to $k$ singular values: Let $A \in \mathbb{R}^{m\times n}$ and assume w.l.o.g. $m \geq n$. Then there exist orthonormal matrices $U \in \mathbb{R}^{m\times m}$ and $V \in \mathbb{R}^{n\times n}$ as well as a diagonal matrix $S = \operatorname{diag}(s_1, ..., s_n) \in \mathbb{R}^{m\times n}$ with singular values $s_1 \geq s_2 \geq \ldots \geq s_n \geq 0$, such thatn

$$A = USV^\top. \tag{14}$$

The columns of $U$ and $V$ are called left and right singular vectors, respectively. The best rank-$k$ representation of the matrix is then given by

$$A \approx \tilde{A} = \tilde{U}\tilde{S}\tilde{V}^\top = \begin{pmatrix} U_1 & \ldots & U_k \end{pmatrix} \operatorname{diag}(s_1, ..., s_k) \begin{pmatrix} V_1 & \ldots & V_k \end{pmatrix}^\top, \tag{15}$$

where $\tilde{U} \in \mathbb{R}^{m\times k}, \tilde{V} \in \mathbb{R}^{n\times k}$ and $\tilde{S} \in \mathbb{R}^{k\times k}$. If we only store the result $\tilde{U}_s = \tilde{U}\tilde{S}$ and $\tilde{V}$, we reduced the required memory from $\mathcal{O}(mn)$ to $\mathcal{O}((m+n)\,k)$.

In this paper, we will concentrate on a pure SVD-based compression by considering the phase-space distribution as a matrix where the rows represent a point of the (up to three dimensional) spacial grid and the columns represent a point of the (up to three dimensional) velocity grid. There exist, however, a number of low-rank tensor formats that generalize the concept of the singular-value decomposition for the compression of high-dimensional data.[73],[74],[75] This can be exploited to further compress also within configuration and velocity space, respectively. We leave the question of which representation of the distribution function yields the best compromise between memory and computational complexity to future work.

### 3.2 Efficient reconstruction of the solution at a given point in time: randomized SVD

In our case we want to avoid storing the full (high-dimensional) distribution function and instead we need an algorithm, which allows for *lazy evaluation* of the full tensor to compute a low-rank representation of it. There are a number of such approaches in the literature and we focus on the *randomized singular value decomposition* (RSVD) in the paper. For the tensor case, alternative algorithmus as the blackbox approximation[76],[77] or the TT cross approximation[78] can also be considered.

RSVD provides an efficient alternative to a full SVD, which also supports lazy evaluation, by using random sampling to identify a subspace that captures the dominant action of the matrix.[79]

The key idea is to multiply the target matrix $A \in \mathbb{R}^{m\times n}$ with a random test matrix $\Omega \in \mathbb{R}^{n\times k}$, producing a sketch $Y = A\Omega$ that spans (with high probability) the dominant column space of $A$. In our experiments we will use a Gaussian random matrix. After orthonormalizing the columns of $Y$ to obtain $Q$, one computes a reduced SVD of the much smaller matrix $B = Q^\top A$. The approximation

$$A \approx QB \tag{16}$$

then provides a rank-$k$ factorization. The procedure is summarized in Algorithm 2.[79]

**Algorithm 2** Randomized SVD with on-the-fly evaluation (transpose formulation)

**Require:** Matrix $A \in \mathbb{R}^{m\times n}$, target rank $k$, oversampling $p$.
1: Draw random Gaussian matrix $\Omega \in \mathbb{R}^{n\times(k+p)}$.
2: Form sketch $Y = A\Omega$ using only matrix-vector products.
3: Orthonormalize: $Y = QR$.
4: Build $B^\top \in \mathbb{R}^{n\times(k+p)}$ column by column:
5: **for** $j = 1, \dots, (k+p)$ **do**
6: $\quad B^\top[:, j] \leftarrow A^\top Q[:, j]$.
7: **end for**
8: Compute SVD: $B^\top = V\Sigma\tilde{U}^\top$.
9: Set $U = Q\tilde{U}$.
10: **return** $U, \Sigma, V$ with $A \approx U\Sigma V^\top$.

In our case the matrix $A$ will be a matrization of the higher dimensional tensor representing the distribution function $f$ on a grid. We can avoid explicitly storing it and use lazy evaluation via providing a function, which evaluates $Ax$ and $A^\top y$ on-the-fly for any $x \in \mathbb{R}^n$ and $y \in \mathbb{R}^m$ at the cost of having to evaluate the same value of $A$ (and thus $f$) several times. Oversampling with $p \ll k$ increases the robustness of the approach. In our simulation a good choice was $3 \le p \le 10$.

## 4 Numerical Study

In this section we want to investigate how the introduction of restarts and low-rank compression influences the accuracy and conservation properties of NuFI. Note that when restarting NuFI via storing the distribution function periodically in time, i.e., every $n_t^r \in \mathbb{N}$ time steps, this makes NuFI comparable to a backwards Semi-Lagrangian (BSL) approach with subcycling. In fact, in the limit $n_t^r \to 1$ this approach becomes just a normal BSL scheme. Thus in addition to investigating how NuFI's accuracy is affected by the restarts, we also will compare NuFI-LR to both a linear and cubic BSL. The former is meant as a direct comparison using the same interpolation order and the later as cubic BSL as comparison to a state-of-the-art approach.

## 4.1 Two stream instability 1D

In previous works on DLR it has been shown that it is able to reproduce kinetic effects such as *Landau Damping* quite well,[56] which however, is not very surprising as the distribution functions associated to Landau Damping remain low rank. A far more challenging kinetic phenomenon to simulate – in the absence of boundary effects and without going to electro-magnetic Vlasov-Maxwell system – is the *Two Stream Instability*, where two counter-streaming beams collide and start mixing. In the simplified one-dimensional space and velocity setting this instability can be for example simulated by initializing the distribution function as

$$f_0(x,v) = \frac{1}{\sqrt{2\pi}}(1+\alpha\cos(kx))v^2\exp(-v^2/2), \tag{17}$$

where $\alpha$ is the strength and $k$ the wave length of the perturbation. We choose $\alpha = 0.01$ and $k = 0.5$ in this case. The boundary conditions are chosen to be periodic with $x \in [0, 4\pi]$.

Even in this simplified setting correctly capturing the evolution of the Two Stream Instability requires a very high phase-space resolution for Semi-Lagrangian schemes: After an initial damping of the perturbation on the beams leading to a slight damping of the electric energy, the beams start interacting and thereby trigger an instability. This leads to an intermediate stage with linear growth of the electric energy, until a saturation point is reached at which point a phase-space vortex has formed. After the saturation point is reached, the electric energy starts oscillating around a fixed level and – due to the lack of collisions in the Vlasov– Poisson system – filamentation starts forming in the distribution function. Especially grid-based methods have a hard time capturing all of the above mentioned effects. Numerical dissipation due to the finite grid resolution leads to inaccuracies in the prediction of the growth rate, but more gravely also prevents correct resolution of the phase-space filamentation for long simulation periods, which in turn leads to a loss of electric energy and unphysically smooth distribution functions. This issue, while less problematic in a synthetic simulation like here, may well lead to unphysical results in long, production simulations, where several kinetic effects overlap and influence each other, and thus it becomes important to resolve them accurately over long time scales.

All simulations – unless indicated otherwise – are carried out on a phase-space grid with $256 \times 256$ degrees of freedom and with a time-step of $\Delta t = \frac{1}{10}$. The cut-off in the velocity space is chosen as $v_{\max} = 6$. As our current prototype-implementation of NuFI-LR employs linear interpolation in between grid-nodes, we decided to also compare to a 1st order in phase-space (backwards) Semi-Lagrangian scheme with 2nd order in time integration[2], which we refer to as *SL* in the following. Moreover, we also compare to a SL scheme with cubic spline interpolation which is considered state-of-the-art for Semi-Lagrangian Vlasov simulations.

First we compare NuFI with NuFI-LR while varying the maximum allowed rank. In our implementation we truncate by prescribing a maximum rank as well as (relative) cut-off tolerance. In the following, if a tolerance is not explicitly mentioned, we use machine precision, i. e., always allow the maximum cut-off rank.

Before looking at the solution quality, we consider the evolution of the singular value spectrum of the phase-space distribution function $f$ in figure 1. In the linear phase $t \leq 10$ the spectrum is quite small with only the first few singular values being relevant, but when the dynamics reach the non-linear phase the decay becomes increasingly shallow. Thus with increasing time the rank of the solution also grows, however, the steepest rank growth is in the instability growth phase between $t \approx 10$ and 25.

[2]In fact we use the reduction of the restarted NuFI scheme to $n_t^r = 1$ as it is equivalent to a 1st order BSL scheme.

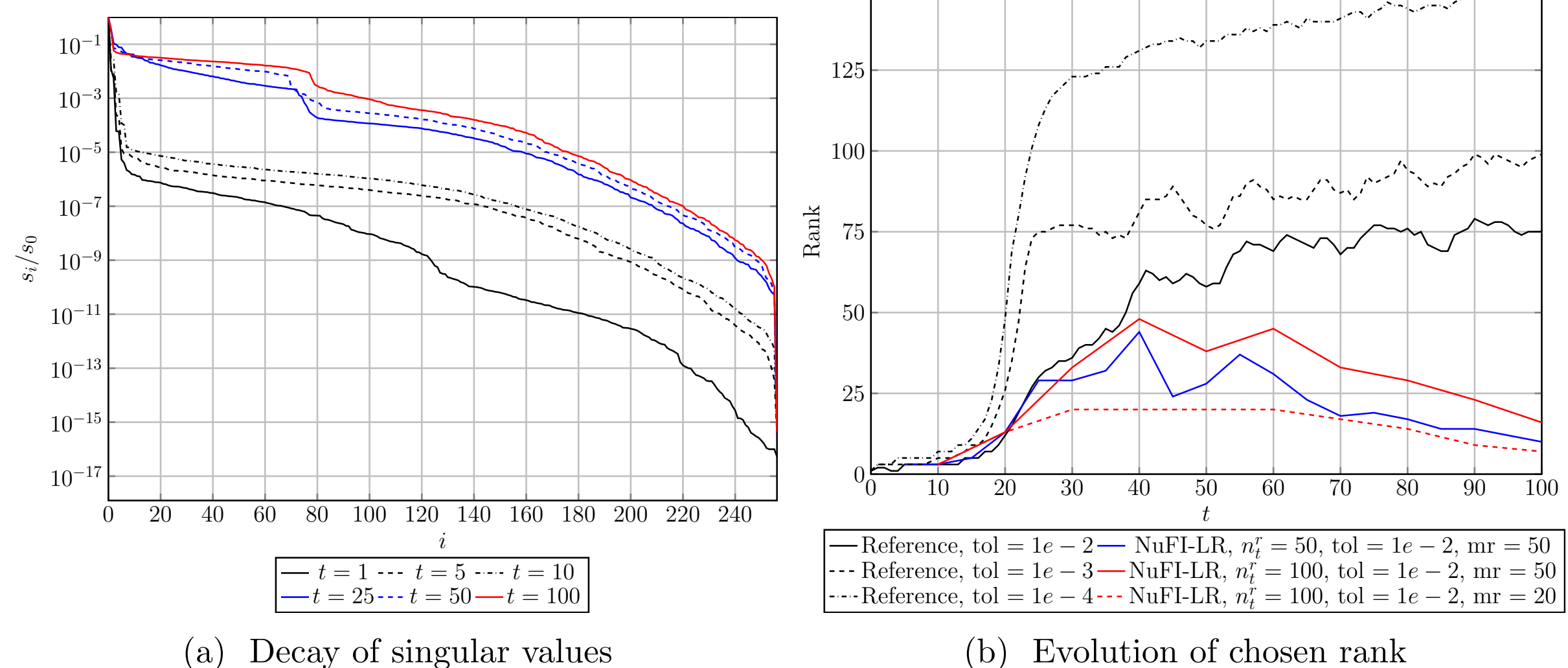

(a) Decay of singular values

(b) Evolution of chosen rank

Figure 1: In this figure we look at the evolution of singular values and ranks in a simulation of the 1-dimensional *two stream instability*. On the left we look at the decay of the (normalized) singular values for a reference simulation with NuFI. Prior to the non-linear stage ($t \leq 10$) the number of relevant singular values is small (black lines). When arriving the at the non-linear stage the decay of the singular values becomes substantially slower (blue and red), however, there is a clear jump between 60-80 depending on the time, meaning that most information is still carried within these first 60-80 modes. The right figure looks at the evolution of the chosen ranks for NuFI-LR (blue/ red) simulations and compares those to cut-offs with (relative) tolerances set to $10^{-2}, 10^{-3}$ and $10^{-4}$. There is a clear increase of required rank which can be seen in the growth phase of the instability between $t = 10$ and 25. In general the chosen ranks for NuFI-LR are at least a factor 2 smaller than the objectively correct choice, however, due to the restarts some information already got lost at that point.

In figure 2 we look at how truncating the maximum rank of the distribution function influences the accuracy of reproducing the evolution of the electric energy. In figure 2a we see that the overall evolution is captured well by all simulations. Even when reducing the maximum rank (mr) to 5 the evolution, including the growth rate of the instability, is reproduced correctly until $t \approx 30$. When zooming into the nonlinear phase after the instability is saturated, we see that indeed the simulations only start deviating after $t \approx 30$. As to be expected, the deviation becomes larger when decreasing the rank, however, all simulations are able to keep oscillating around the correct energy level even if the simulations with the lowest rank of 5 and 10 almost "flatline" for $t \to 100$. The period of the large oscillation seems to be captured correctly by all simulations but the amplitude decreases with decreasing rank.

Next we keep the maximum allowed rank fixed to 20 and vary the restart period, see figure 3. We observe that with increasing $n_t^r$ the accuracy of the simulation also increases. For $n_t^r = 100$ and 50 the results are similar and remain close to the reference solution, while decreasing $n_t^r$ to 10 or 5 introduces significant numerical dissipation, which manifests itself through strong loss in amplitude and induces an unphysical drift in the electric energy. For small values of $n_t^r$ the growth rate of the instability is slightly smaller. Additionally, when comparing the figures 3b and 3d, we see that in fact the frequency of restarts seems to have a greater impact on physical accuracy: choosing high rank but a high restart frequency, i.e., small $n_t^r$, leads to worse results than choosing low rank with a low restart frequency. This suggests that high compression rates can be achieved while retaining good physical accuracy with NuFI-LR – comparable to full rank NuFI solutions – as long as $n_t^r$ is large enough.

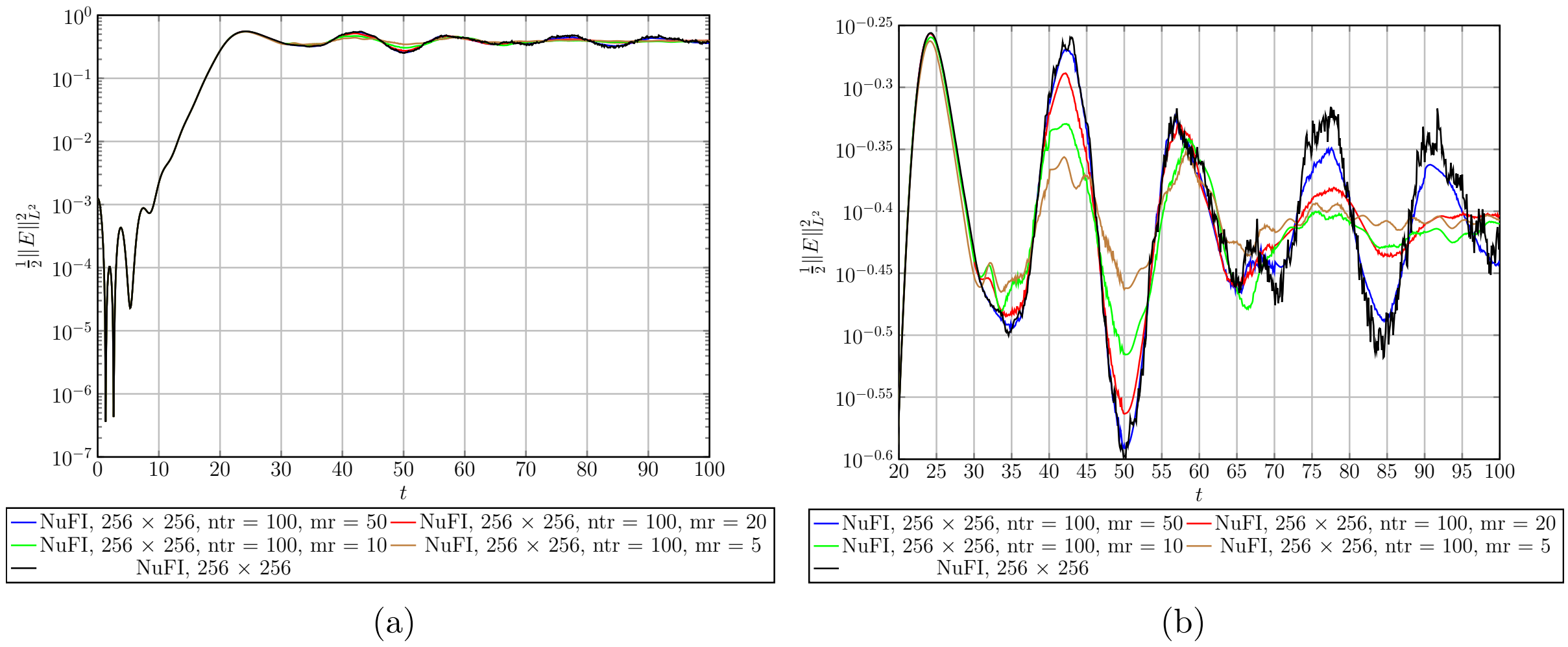

Figure 2: In this figure we compare the evolution of the electric field for the *two stream instability* in $d = 1$ between NuFI and NuFI-LR. The restart period is set to $n_t^r = 100$ for all NuFI-LR simulations and we vary the maximum rank (mr) of the solution. In the left figure we show the full evolution until $t = 100$, while on the right we zoom into the non-linear mixing phase after $t \approx 20$. NuFI and NuFI-LR agree very well until about $t \approx 30$ and, in particular, all simulations are able to capture the correct growth rate of the instability. In the later stage of the simulation, when decreasing the maximum rank of the solution, the amplitude of the oscillation drops and the very low rank solutions exhibit some noise, however, all simulations are able to capture the correct energy level until $t = 100$.

In figure 4 we compare the evolution of the electric energy between NuFI, NuFI-LR and SL. Similar to what has been already observed in previous work[69] one requires a higher resolution to achieve the same accuracy with SL using the same interpolation order as for NuFI. In this case linear SL requires a $1024 \times 1024$ grid, i.e., a 4 times higher phase space resolution, to achieve the same accuracy in predicting the growth rate as NuFI with $256 \times 256$ resolution in phase space. Still even in that case the high resolution linear SL solution starts to loose amplitude faster than NuFI and already exhibits a small unphysical drift in the electric energy in the late stage of the simulation. When choosing the same resolution for linear SL as for NuFI, linear SL predicts a smaller growth rate and exhibits the energy drift directly from the start of the non-linear stage after $t \approx 25$. At $t = 100$ linear SL dissipated about 75% of the electric energy in the system. Compared to this the worst energy drift observed for NuFI-LR is about 25% for a compression down to max rank of 20 and $n_t^r = 5$, see figure 3b, which is at a similar level as the drift for the high resolution linear SL. Note that the cubic SL simulation remains more accurate than its linear counterpart preserving the energy level until the end of the simulation. The accuracy of the cubic SL is comparable to the full-rank NuFI and highest rank NuFI-LR simulations with less noise but a subtle loss in amplitude.

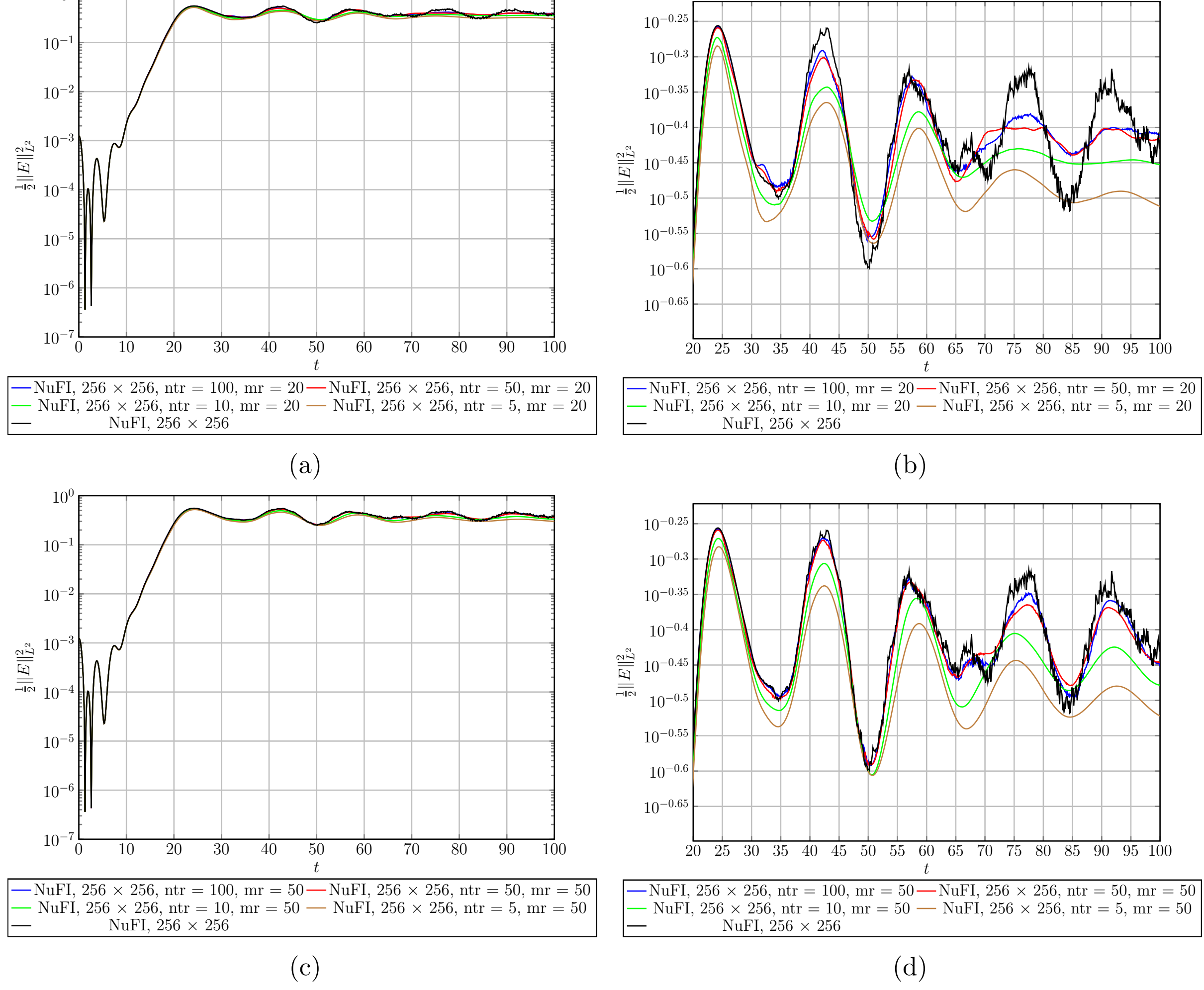

Figure 3: In this figure we compare the evolution of the electric field for the *two stream instability* in $d = 1$ NuFI with NuFI-LR for varying restart frequencies. The figures in the top row show results with a maximum rank of 20 and in the bottom with a maximum rank of 50. As in the last figure we show the full evolution on the left and zoom into the non-linear phase on the right. As before we see that the choice of maximum rank mostly affects the amplitude and level of noise in the solution. A higher restart frequency leads to an unphysical drift in the electric energy level.

As expected, the compressed version of NuFI, NuFI-LR, shows a loss of amplitude compared to the full rank solution. NuFI-LR also exhibits noise in form of additional, unphysical oscillations. However, even the lowest rank NuFI-LR simulation reproduces the evolution of electric energy substantially better than the full SL simulation with same interpolation order and it's accuracy is comparable to the cubic SL simulation.

In figure 5 we compare the distribution functions at $t = 50$ between SL, NuFI and NuFI-LR with max ranks of 20 and 50. While SL at this point has already dissipated away almost all filamentation in the distribution function and only a vortex is still visible, for NuFI the filamentation remains very clearly visible. The low rank compression in NuFI-LR also leads to some loss of structure and noise, but to a lesser degree than for SL. In fact even for compression down to max rank 20 most of the filamentation remains visible, albeit less clearly than for the full rank solution. With max rank set to 50 the filamentation in the core of the vortex is still clearly visible and only in the boundary region the filamentation is slightly smeared out. The cubic SL simulation is able to preserve far more structure than its linear counterpart, however, still to a lesser degree than both the full rank NuFI and mr = 50 simulation with NuFI-LR.

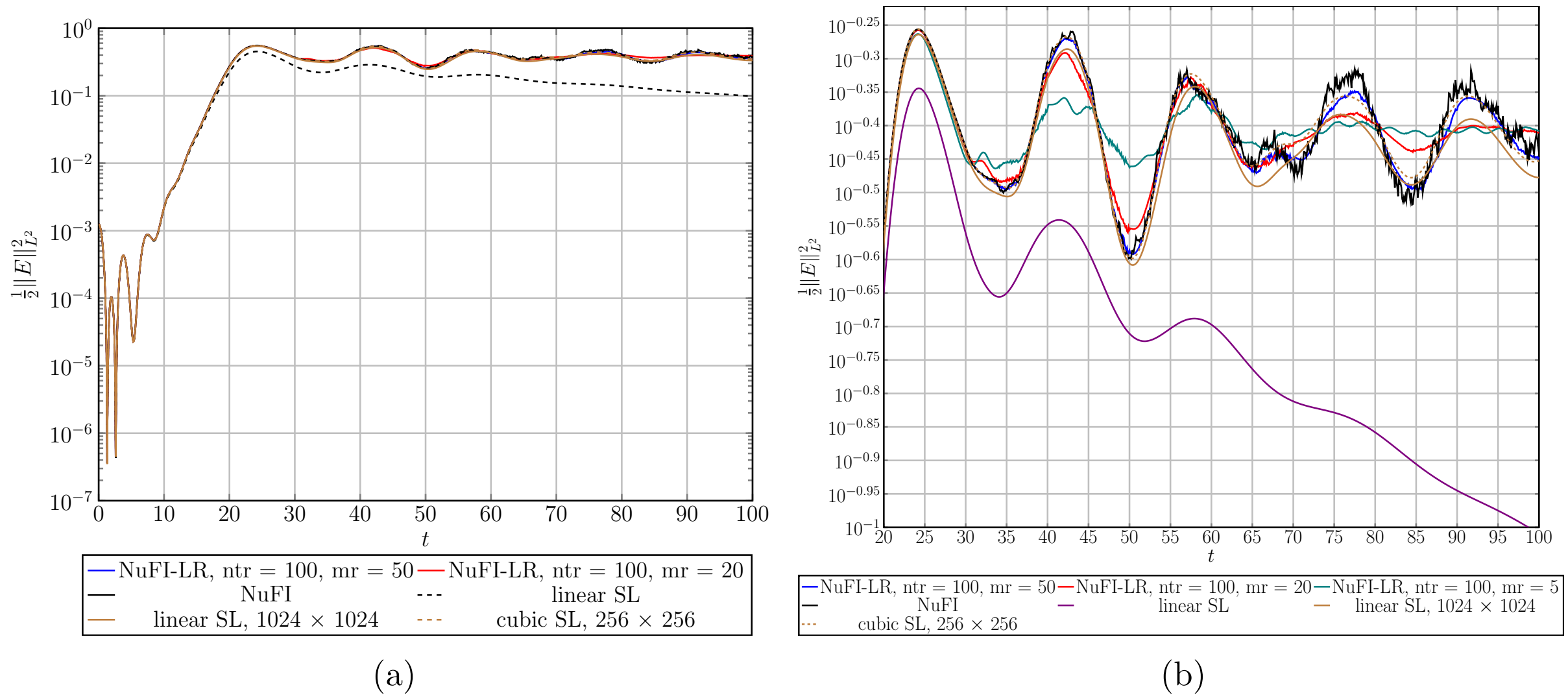

Figure 4: In this figure we compare the evolution of the electric field for the *two stream instability* in $d = 1$ between NuFI-LR and Semi-Lagrangian (SL) and Semi- Lagrangian with low rank compression (SL-LR). The comparison includes NuFI-LR with $n_t^r = 100$ and the compression is performed to a maximum ranks of $5, 20$ and $50$. To verify the results we also included a high resolution SL simulation. The left figure shows the full evolution and the right figure zooms into the non-linear phase after the instability is saturated. NuFI-LR agrees with the reference solution until $t \approx 30$, even with a maximum rank as low as 5, while already the full rank SL solution starts deviating earlier than $t = 20$ and fails to capture the correct growth rate. When introducing low rank compression for SL, the error is substantially larger for the very low rank solution with maximum rank of 5. The noise is much lower in SL for the higher rank solutions than for NuFI(-LR), which is due to the numerical diffusion caused by the underlying grid and also causes a strong drift in the electric energy. Choosing the resolution for SL 4 times higher in each direction allows to reproduce NuFI's full rank result, but is still less accurate in the prediction of the electry energy level for the late stages of the simulation.

Finally, we also consider the conservation of (total) Energy and Entropy in figure 6. Recall that the pure NuFI algorithm preserves entropy as well as $L^p$-norms analytically, i. e., to machine precision and total energy up to the discretization error but without drift.[66] We observe that compression introduces a loss of conservation properties in NuFI-LR, however, the error grows depends on the restart frequency, while the compression rate itself – for the chosen ranks – seems to have only a relatively small influence. When choosing $n_t^r = 100$ the error made by NuFI-LR in the conservation of total energy is less than 2% until $t = 1000$. For $n_t^r = 10$ it goes up to 12%. For entropy the error stays below 4% for $n_t^r = 100$ and goes up to 13% for $n_t^r = 10$. Meanwhile for SL the errors in total energy conservation saturate at 29% and for entropy at around 22%.

Finally, let us note that a low-rank compression of the Semi-Lagrangian scheme would also be possible to speed up the simulations, e.g. through DLR. However, those results cannot be more accurate than the underlying Semi-Lagrangian scheme and thus suffer from the same reduced accuracy as the pure SL scheme compared to NuFI-LR.

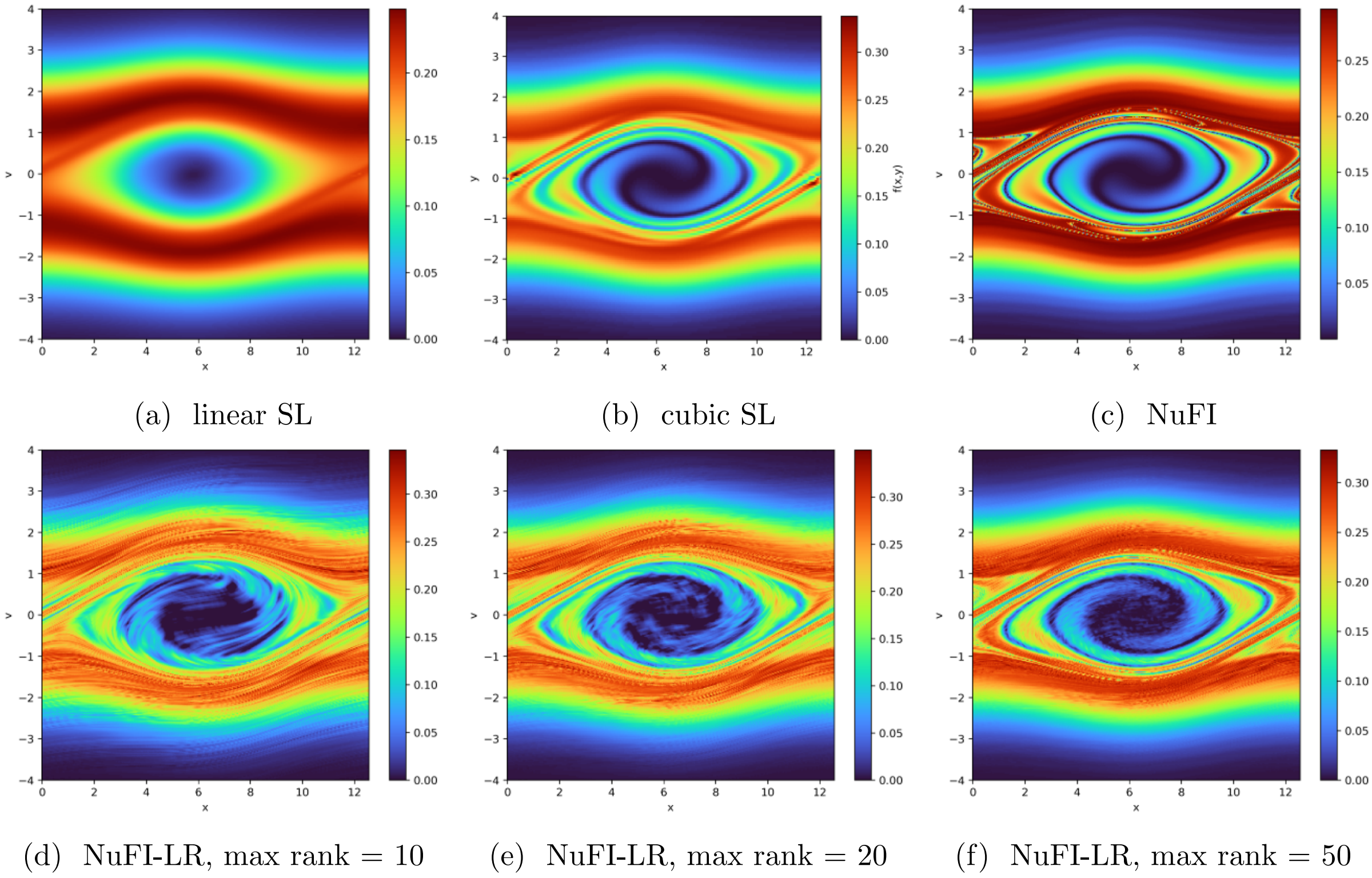

(a) linear SL (b) cubic SL (c) NuFI

(d) NuFI-LR, max rank = 10 (e) NuFI-LR, max rank = 20 (f) NuFI-LR, max rank = 50

Figure 5: Here we compare the distribution functions $f$ at time = 50 of simulations of the *two stream instability* in $d = 1$ performed by (linear/cubic) SL, NuFI and NuFI-LR with the same base resolution of $256 \times 256$ in phase-space. At this point the numerical diffusion in linear SL (top left) already smeared out almost all details in $f$ and only the vortex remains visible, however, also with reduced maximum height of the solution. The cubic SL (top middle) is able to preserve more details but still less clearly than the full rank NuFI simulation (top right), which is able to reproduce all the filamentation present in the distribution function. The high rank NuFI-LR (mr = 50, bottom right) is also able to capture more details than cubic SL but less clearly than full rank NuFI. With moderate rank (mr = 20, bottom middle) NuFI-LR is capable to reproduce most of the larger filaments, though fails at the boundary and shows some numerical artefacts, comparable to the cubic SL results. The lowest rank NuFI-LR (mr = 10, bottom left) shows a lot of numerical artifacts with quality comparable to the full rank linear SL.

## 4.2 Ion-Acoustic shock with reflecting wall

One of the major advantages of NuFI over grid-based approaches is that handling of non-periodic boundary conditions works similar to particle-based schemes and therefore is simpler to implement. For DLR, while possible, non-periodic boundary conditions introduce additional complications.[80] In this section, we want to showcase on a simple example of an ion-acoustic shock induced by a reflecting wall that NuFI-LR is capable to deal with non-periodic boundaries as well. The setup is taken from previous work of one of the authors[69] and is built upon a simplified version of shock simulation suggested by Liseykina et al.[81] We again restrict ourselves to $d = 1$ with only electro-static forces. This time we consider both electrons and ions (Hydrogen-ions), where we set the realistic mass ratio of $M_r = 1836$, i. e., $m_e = 1$ and $m_i = 1836$. The physical domain is from $x_{\min} = -L$ to $x_{\max} = 0$ with $L \gg 1$. The ions and electrons are initialized with unperturbed Maxwellian velocity distributions, i. e., the initial distribution for $\alpha \in \{i, e\}$ are

$$f^\alpha(x, v) = \sqrt{\frac{m^\alpha}{2\pi T^\alpha}} \exp\left(-\frac{m^\alpha(v-u_s^\alpha)^2}{2T^\alpha}\right), \tag{18}$$

where the ion distribution is centered around a drift speed $u_s > 0$. The right boundary is set to be a perfectly reflecting wall for both electrons and ions. The left boundary is an open boundary for particles leaving the domain and for particles entering the domain prescribes an in-flow with the same Maxwellian as the initial condition for the particle species. Both boundaries are set to be perfectly conducting, i. e., for the Poisson equation we set zero Neumann-boundary conditions. Following the suggestions of Liseykina et al. in the $u_s = 0.4$ case we choose $T_e = 11475 \gg T_i = 10$ and $L = 200$. We simulate with NuFI choosing $N_x = 1024$, adaptive integration in velocity with at least $N_v = 64$ cells and a time-integration step of $\Delta t = \frac{1}{10}$. To restart we use $n_t^r = 50$, and $N_x^r = N_v^r = 1024$. Additionally we compress to a maximum rank of 20 and set a relative cut-off tolerance of $10^{-3}$ for the singular values.

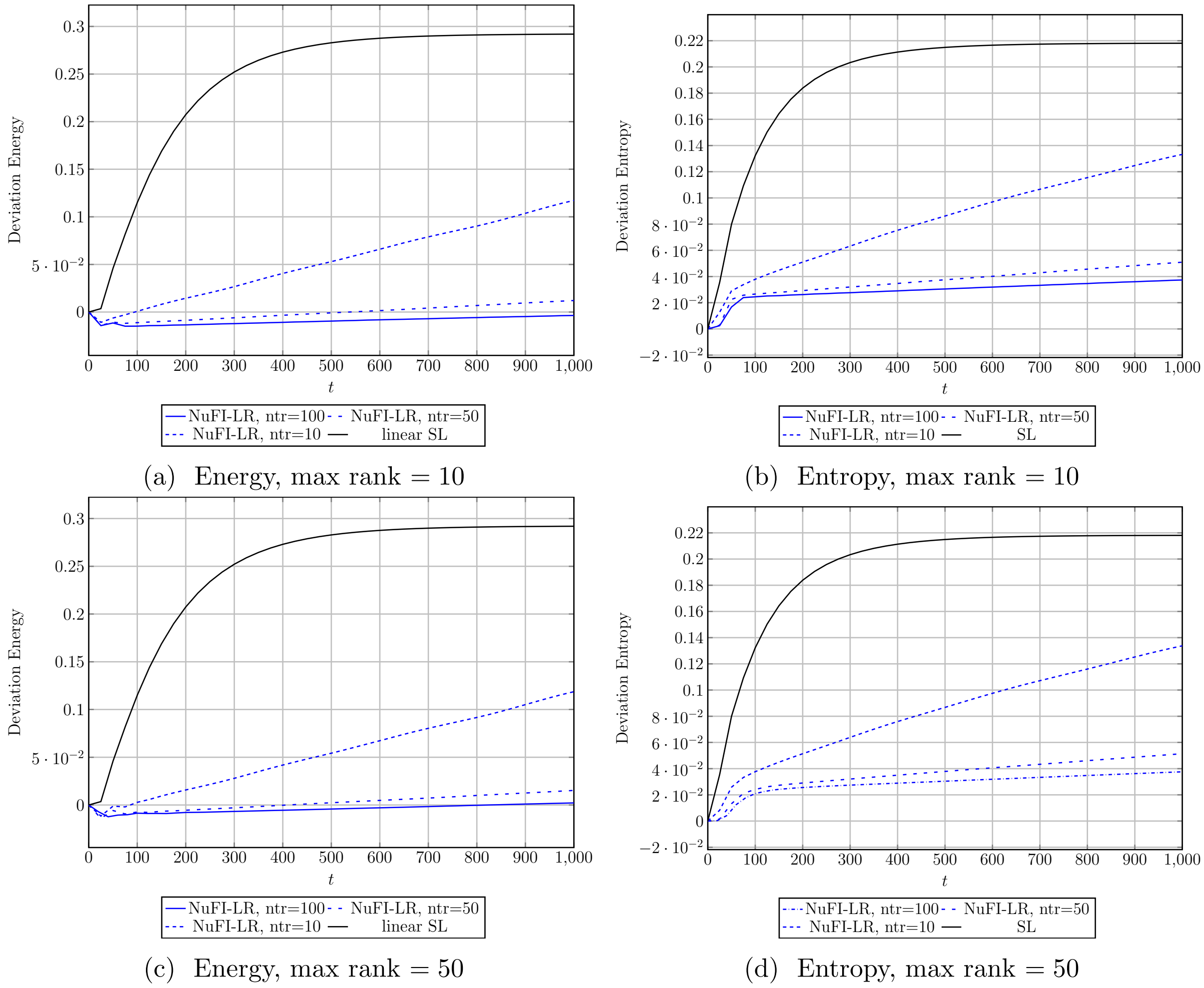

(a) Energy, max rank = 10 (b) Entropy, max rank = 10

(c) Energy, max rank = 50 (d) Entropy, max rank = 50

Figure 6: In this figure we compare the conservation of total energy (left) and entropy (right) between simulations of the one-dimensional *two stream instability* via compressed versions of NuFI and SL. Overall we see that NuFI-LR conserves these quantities much better for all compression rates and that the conservation error depends almost exclusively on the restart frequency in this test case. The conservation error at $t = 1000$ is about an order of magnitude larger when decreasing $n_t^r$ from 100 to 10.

The resulting electron and ion distribution function at $t = 100, 500$ and $1000$ are shown in figure 7. As we already observed in the previous section, NuFI-LR is able to capture the initial shock formation as well as to preserve the fine structure in the distribution functions, especially the late stage electron distribution function in figure 7b and 7c, to a high accuracy. Additionally, NuFI-LR is also compatible with adaptive integration in velocity space as well as tracking of changes in the velocity support of the distribution functions.[67] Note that over the course of the simulation the velocity support of the distribution functions, especially the ions, grows by more than an order of magnitude. Tracking this phenomena dynamically with SL is hard to implement.

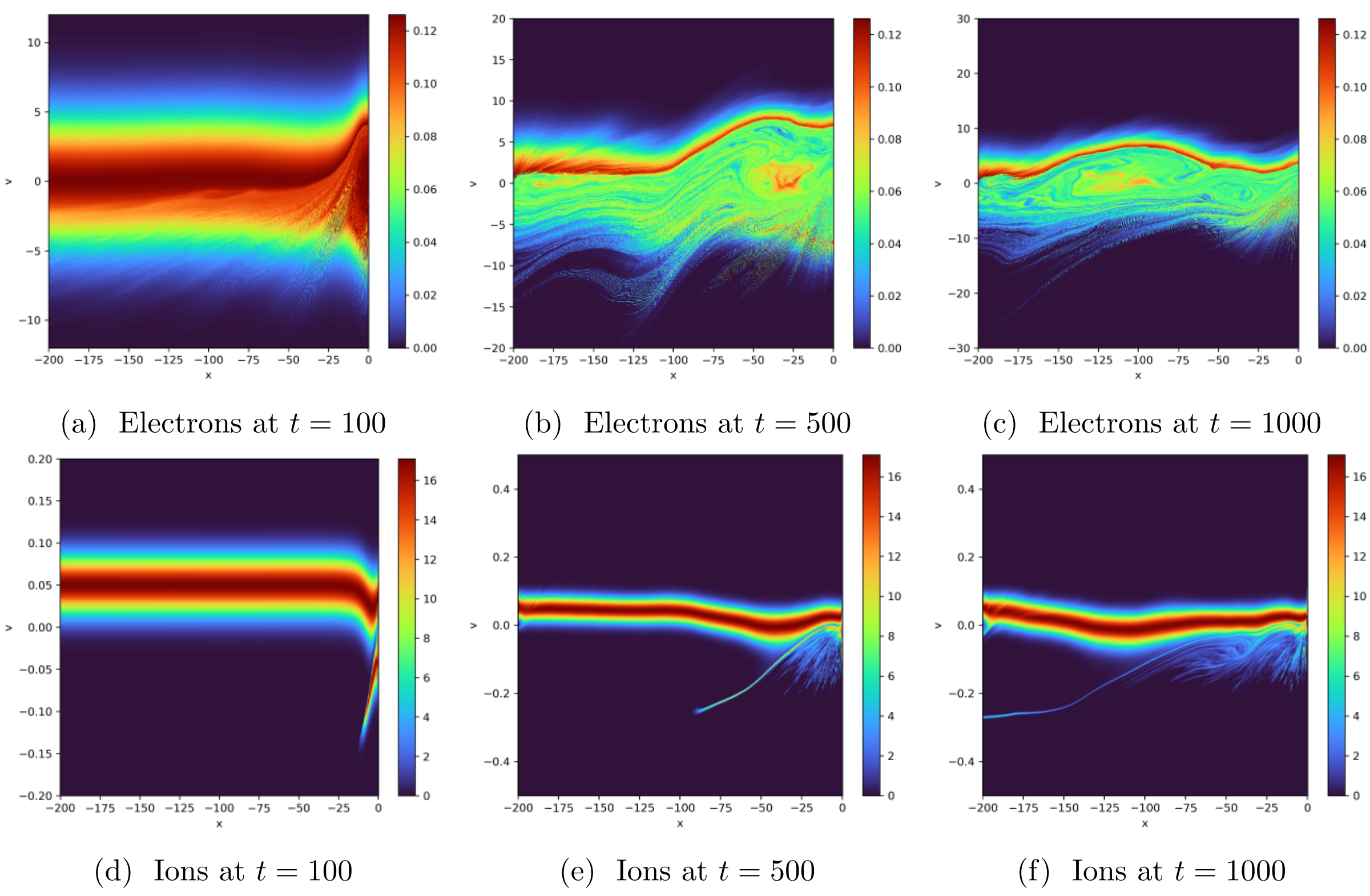

(a) Electrons at $t = 100$ (b) Electrons at $t = 500$ (c) Electrons at $t = 1000$

(d) Ions at $t = 100$ (e) Ions at $t = 500$ (f) Ions at $t = 1000$

Figure 7: The electron (top) and ion (bottom) distribution functions from a simulation of an *ion-acoustic shock* caused by introducing a reflective wall at $x = 0$. The simulation was run with NuFI-LR with maximum rank of 20. No smearing or obvious numerical artifacts appear appear in the simulation time and the fine structures forming around the shock seem to be well resolved. Note that the counter-streaming waves at the left boundary are caused by the finite domain size in $x$ direction and not an effect caused by the compression.

### 4.3 Two stream instability 2D

Finally we also look at a simulation with 4-dimensional phase-space, i. e., 2 space and 2 velocity dimensions. For the same reasons as listed in section 4.1, we decided to again look at a variant of the *two stream instability.* We follow the setup proposed by in previous work of one of the authors.[56] The initial condition is set to

$$f_0(x,v) = \frac{1}{8\pi}\left(1+\alpha\left(\cos(kx)+\cos(ky)\right)\right)\left(\exp\left(-\tfrac{(u-v_0)^2}{2}\right)+\exp\left(-\tfrac{(u+v_0)^2}{2}\right)\right) \\ \left(\exp\left(-\tfrac{(v-v_0)^2}{2}\right)+\exp\left(-\tfrac{(v+v_0)^2}{2}\right)\right), \tag{19}$$

where the perturbation strength is $\alpha = 10^{-3}$, the wave number of the perturbation is $k = 0.2$ and the drift velocity of the beams is $v_0 = 2.4$. The resolution was set to $N_x = N_y = N_u = N_v = 128$ and the time-step to $\Delta t = \frac{1}{10}$.

Similar to the 1d case, the agreement between NuFI and NuFI-LR are overall good, especially in the initial phase until $t \approx 30$ the results are practically identical to the reference solution. The simulations with smaller $n_t^r$ show a larger deviation and a slight drift after $t \approx 32$. Also – in contrast to the 1d case – here we see that choosing a smaller maximum rank led to an earlier deviation from the reference result and led to additional unphysical "jumps" at $t \approx 40$ (restart timing) and $t = 44.5$. However, the results by NuFI-LR – independent of the compression rate and restart frequency chosen here – are still much closer to the reference solution than a full rank SL solution with same resolution. SL already fails to capture the correct growth rate and starts deviating as early as $t \approx 15$.

**Remark 4.1** *The compression for this simulation was performed using RSVD on a matrization of the data tensor. The degrees of freedom are combined into the rows for the spatial and into the columns for the velocity variables, which neglects potential additional structure in the combined directions. However, we chose RSVD over HT compression in this case due to the significantly faster evaluation, even though HT potentially can provide a more memory-slim decomposition. Hence a row/column max rank of* $n$ *implies that the effective rank per direction* $m$ *is likely much smaller, i. e.,* $m \ll n$. *This is an effect present in the 4-dimensional but not 2-dimensional phase-space.*

## 5 Discussion and outlook

In this study we investigated the suitability of low-rank tensor approximation for compression of snapshots of the Vlasov equation to be used as a restart mechanism for NuFI to reduce the computational complexity from quadratic back to linear in the number of time steps. We focused on the example of electro-static *two stream instabilities* in 2- and 4-dimensional phase-space to test how well NuFI-LR can deal with filamentation forming in the distribution function, which is a common phenomenon in both electro-static and fully electro-magnetic Vlasov simulations and one of the major challenges for numerical solvers. Additionally, we demonstrated on the example of an ion-acoustic shock conditions that NuFI-LR is also applicable to more complex setups such as non-periodic boundary conditions and remains feasible as well as stable for longer simulation periods.

The results from sections 4.1 and 4.3 clearly show – for the presented test cases – that while low-rank compression reduces the accuracy of the NuFI approach, the algorithm remains substantially more accurate than a Semi-Lagrangian approach using comparable interpolation order and resolution. In our numerical experiments, the restart interval has been more decisive for the quality of the solution than the compression rate. In particular, in the nonlinear phase which is characterized by micro-scale filamentation, the restart interval supposedly needs to be large enough to resolve the time scale of micro-scale processes of the physical system under consideration while the exact detail of the distribution functions at the restart distribution seems to be less important to recover the statistics correctly. Because NuFI uses symplectic integration in time it preserves the Hamiltonian structure of the Vlasov solution, which is broken when restarting. Hence frequent restarts lead to a build-up of the associated error. Additionally, projecting on a phase-space grid (with linear interpolation) introduces strong numerical diffusion in the system.

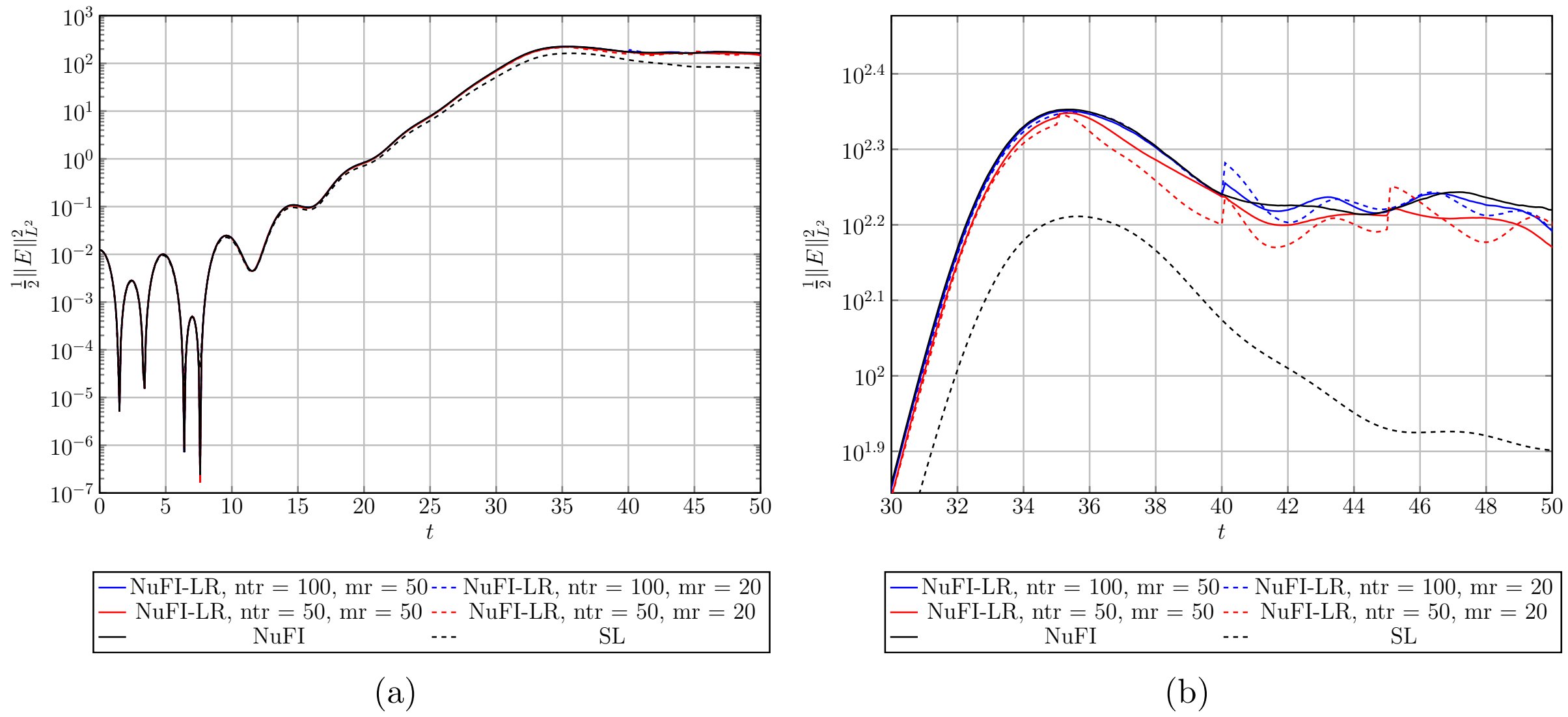

Figure 8: In this figure we look the evolution of the electric energy for simulations of a *two stream instability* in $d = 2$. We again compare the results between NuFI, NuFI-LR and SL. As alreay observed in the one-dimensional case NuFI and NuFI-LR both are more accurate than SL in predicting the evolution. However, in contrast to the one-dimensional case we see more noise in the max rank 20 solution which is likely due to the compression being computed through a RSVD on a matrization of the 4D data, which depend on the choice of matrization may miss structure and also reduces the effective rank per dimension.

The choice of compression algorithm in our study was very rudimentary and only meant as a proof-of-concept. To find the best suited candidate one has to perform a study of advanced state-of-the-art compression algorithms in terms of accuracy and speed. In particular, this study should also include data from fully electro-magnetic simulations and, if possible, also observational data to verify the suitability of the approach in all possible cases.

Future work will compare the use of various tensor formats to investigate which format is most efficient for the purpose of restarting NuFI. In particular, for our purpose the requirement on the low-rank compression is very specific: not only the compression step has to be fast and accurate, while allowing for lazy evaluation, but also the compressed data format has to allow for fast point-wise evaluation of the data or at least fast evaluation of small data slices. Commonly algorithms and formats in literature are optimized towards fast compression and decompression of large data-subsets, while point-wise evaluation may not be required. This algorithmic challenge or change of perspective has to be tackled to make NuFI-LR a viable alternative for production simulations.

## Acknowledgments

R.P.W. wants to thank the organizers of the National High Performance Computing (NHR) graduate school and the Federal Ministry of Education and Research as well as the state governments for supporting this work as part of the joint funding of National High Performance Computing (NHR). R.P.W. wants to thank the DFG, who funded part of the work through the IRTG Modern Inverse Problems (333849990/GRK2379). This work has received funding from the European High Performance Computing Joint Undertaking (JU) and Belgium, Czech Republic, France, Germany, Greece, Italy, Norway, and Spain under grant agreement No 101093441. Views and opinions expressed are however those of the author(s) only and do not necessarily reflect those of the European Union or the European High Performance Computing Joint Undertaking (JU) and Belgium, Czech Republic, France, Germany, Greece, Italy, Norway, and Spain. Neither the European Union nor the granting authority can be held responsible for them. K.K. acknowledges support by the Deutsche Forschungsgemeinschaft (DFG, German Research Foundation)—project ID 530709913.